\documentclass{amsart}
\usepackage[utf8]{inputenc}
\usepackage{amsmath}
\usepackage{amssymb}
\usepackage{amsthm}
\usepackage{tkz-graph}
\usepackage[numbers]{natbib}
\usepackage{enumerate} 
\usepackage{bm} 
\usepackage{varwidth}
\usepackage{listings}
\usepackage{graphicx}
\usepackage{color}
\usepackage{pstricks}
\usepackage{calc}
\usepackage{float}
\usepackage{hyperref}
\theoremstyle{definition}

\newcommand{\NN}{\mathbb{N}}

\newcommand{\Cr}{\operatorname{cr}}

\newcommand{\outdeg}{\operatorname{outdeg}}

\newcommand{\cupdot}{\mathbin{\mathaccent\cdot\cup}}
\newcommand{\apex}{\operatorname{apex}}
\newcommand{\coapex}{\operatorname{coapex}}
\newcommand{\Act}{\operatorname{active}}
\newcommand{\faces}{\operatorname{faces}}
\newcommand{\abs}[1]{\lvert #1 \rvert}
\newcommand{\diredge}{%
\begin{tikzpicture}%
\draw[fill = black] (.25ex,.25ex) circle (.3ex);
\draw[thick, ->] (.55ex,.25ex) -- (1.85ex,.25ex);%

\end{tikzpicture}%
}

\title[Comp. classificaltion of some almost minimal Ramsey graphs]{A computerised classification of some almost minimal triangle-free Ramsey graphs}
\author{Oliver Kr\"uger}
\date{\today}

\begin{document}

\maketitle

\begin{abstract}
A graph $G$ is called a \emph{$(3,j;n)$-minimal Ramsey graph} if it has the least amount of edges, $e(3,j;n)$, given that $G$ is triangle-free, the independence number $\alpha(G) < j$ and that $G$ has $n$ vertices. Triangle-free graphs $G$ with $\alpha(G) < j$ and where $e(G) - e(3,j;n)$ is small are said to be almost minimal Ramsey graphs.

We look at a construction of some almost minimal Ramsey graphs, called $H_{13}$-patterned graphs. We make computer calculations of the number of almost minimal Ramsey triangle-free graphs that are $H_{13}$-patterned. The results of these calculations indicate that many of these graphs are in fact $H_{13}$-patterned. In particular, all but one of the connected $(3,j;n)$-minimal Ramsey graphs for $j \leq 9$ are indeed $H_{13}$-patterned.
\end{abstract}

\section{Introduction}
\label{introduction}

In this report we introduce the notion of a type of small directed graphs, called crochet patterns, due to J. Backelin (see \cite{carc}) from which one can build (using a recursive procedure implemented on a computer) larger graphs with interesting properties. The graphs that we are interested in are graphs that are either minimal Ramsey graphs or close to being so.

Let $G = (V,E)$ be a graph (unless otherwise stated a graph will mean a simple undirected graph). Let $n(G)$ and $e(G)$ denote the number of vertices and edges of $G$, respectively. The \emph{clique number} of $G$, $\omega(G)$, is the maximum size of a complete subgraph of $G$. The \emph{independence number} of $G$, $\alpha(G)$, is the maximum size of an independent set of vertices in $G$.

An \emph{$(i,j; n,e)$-graph} is a graph $G$ such that $\omega(G) <i$, $\alpha(G) < j$, $n(G) = n$ and $e(G) = e$. An $(i,j;n)$-graph is a graph that is an $(i,j; n,e)$-graph for some $e \in \NN$ and an $(i,j)$-graph is a graph that is an $(i,j; n)$-graph for some $n \in \NN$. Ramsey's theorem states that for all $i$,$j$ there is some $n$ such that there are no $(i,j;n)$-graphs. The least such $n$ is the (two-colour) \emph{Ramsey number}, $R(i,j)$. The $(i,j)$-graphs are sometimes called Ramsey graphs.

Let $e(i,j;n)$ denote the minimum number of edges in an $(i,j;n)$-graph. Note that the values $e(i,j;n)$ determines the Ramsey numbers since $R(i,j) = \min\{n \in \NN \mid e(i,j;n) = \infty\}$. These numbers are called the \emph{minimum edge-numbers}. The $(i,j;n,e(i,j;n))$-graphs are called \emph{minimal Ramsey graphs}.

In this report we will focus on the construction of graphs that are either minimal Ramsey graphs or which are $(i,j;n,e)$-graphs for which $e-e(i,j;n)$ is small, where $i = 3$. All graphs we consider will therefore be assumed to be triangle-free.

B. McKay has computed all $(3,4)$-,$(3,5)$- and $(3,6)$-graphs, which are available online at \cite{McKayMinimalRamsey}. We will compare the lists of graphs we compute with these lists in the case when we construct $(3,j)$-graphs for $4 \leq j \leq 6$.

Moreover, J. Goedgebeur and S.P. Radziszowski have computed all the minimal $(3,7)-$, $(3,8)$-, $(3,9)$- and $(3,10)$-graphs and several more (possibly incomplete) lists of $(3,j)$-graphs for $j \geq 11$ (see \cite{goed-radz}).
Their lists of graphs are available online at \cite{HouseofGraphsMinimalRamsey}. Also these lists will be used for comparison with those we compute.

We have summarised some highlights of the result in Table \ref{maintable}. We will see that for some specific values of $j,e$ and $n$ we have that many of the $(3,j;n,e)$-graphs are $H_{13}$-\emph{patterned} (in the sense that they can be obtained by a $H_{13}$-decorated pattern, which we define in the next section).

\begin{center}
\begin{table}
\begin{tabular}{c|c|c|c|c}
$(3,j;n,e)$ & $e - e(3,j;n)$ & \# of $(3,j;n,e)$-graphs & $H_{13}$-patterned & unpatterned \\
\hline
$(3,6;16,32)$ & 0 & 5 & 4 & 1 \\
$(3,7;19,37)$ & 0 & 11 & 11 & 0 \\
$(3,7;20,44)$ & 0 & 15 & 15 & 0 \\
$(3,7;21,51)$ & 0 & 4 & 4 & 0 \\
$(3,8;23,49)$ & 0 & 102 & 102 & 0 \\
$(3,8;24,56)$ & 0 & 51 & 51 & 0 \\
$(3,9;26,54)$ & 2 & 444 & 436 & 8 \\
$(3,9;27,61)$ & 0 & 700 & 700 & 0 \\
$(3,9;28,68)$ & 0 & 126 & 126 & 0 \\
$(3,10;29,59)$ & 1 & 1364 & 1304 & 60 \\
$(3,10;30,66)$ & 0 & 5084 & 5082 & 2 \\
$(3,10;31,73)$ & 0 & 2657 & 2657 & 0 \\
$(3,11;32,64)$ & 1 \footnotemark & $\geq$ 3950\footnotemark & 3403 & $\geq$ 547
\end{tabular}
\caption{Table of some highlights of the counts of $H_{13}$-patterned graphs compared to the known Ramsey-minimal graphs.} \label{maintable}
\end{table}
\end{center}
\addtocounter{footnote}{-1}
\footnotetext{That $e(3,11;32) \geq 63$ was proved by A. Lesser in \cite{lesser}. That this inequality is in fact an equality is shown in \cite{carc}.}
\addtocounter{footnote}{1}
\footnotetext{This is possibly less than the actual number of $(3,11;32,64)$-graphs. The 3950 known ones are either patterned (3403 of these), two disjoint copies 
of $(3,6;16,32)$-graphs, at least one component being the unpatterned one (5 of these), or a union of a $H_{13}$ (1 orbit) with a $(3,7;19,37)$-graph (total of 140 orbits) with an extra edge between the two (140 of these), or $H_{13}+G'$ where $G'$ is one of the $(3,7;19,38)$-graphs that are not $H_{13}$-patterned (402 of these).}

\section{Crochet patterns and patterned graphs}

This description of patterned graphs is as far as we know invented by Backelin, and all the notions in this section is purely due to him and is taken from his work in \cite{carc}. What is presented here is only a particular specialisation of a more general framework. We just give the parts that are needed to explain how the patterned graphs that are enumerated in Table \ref{maintable} and the tables in the appendix were computer-generated.

\subsection{Stitches}

The graphs that we are interested in are such that they may be defined by a sequence of so called stitches from the empty graph. We have several types of stitches, the most fundamental of which is the \emph{$d$-stitch}, which is a form of generalised $d$-extension (in the sense of Radziszowski and Kreher in \cite{radziszowski-kreher91}). By $G_v$ we denote the graph obtained from the graph $G$ by removing the vertex $v \in V(G)$, all neighbours of $v$ and all edges incident to those vertices. A $d$-stitch is a pair $(G,v)$ where $G$ is some graph and $v \in V(G)$, such that $\deg(v) = d$, $\alpha(G_v) = \alpha(G) - 1$ but $\alpha(G-vu) > \alpha(G)$ for all $vu \in E(G)$.

For example, the path with two vertices, $P_2$, is a 1-stitch of the empty graph.

The idea is to study how the graph $G$ can be constructed from the graph $G_v$ in such a way that this is done by adding ``many'' vertices, but only a ``few'' new edges. The procedure described below is an attempt to make this philosophy concrete. 

A subset of vertices $M \subseteq V(G)$ of a graph $G$ is called a \emph{destabiliser} if every maximum independent set of $G$ intersects $M$, i.e. the induced subgraph on $V(G)\setminus M$ has a lower independence number than $G$. A destabiliser is said to be \emph{minimal} if there is no destabiliser properly contained in it.

If $(G,v)$ is a $d$-stitch we call $G$ the \emph{stitch-graph} and $v$ the \emph{apex} of the stitch. Let $N_i(x)$ denote the set of vertices at distance exactly $i$ from $x \in V(G)$. The set $N_2(v) \subseteq V(G)$ is then called the \emph{base} of the stitch.

If $(G,v)$ is a $2$-stitch then $N_2(v)$ induces a bipartite graph (each part being adjacent to different neighbours of $v$) which is a destabiliser of $G_v$. In fact, if $H$ contains a bipartite minimal destabiliser $M\subseteq V(H)$, this set can be used as the base of some 2-stitch, $(G,v)$ where $G_v = H$. In this case we denote the stitch-graph $G$ by $\Cr_{2}(H,M)$.

A graph $H$ is called edge-critical if $\alpha(H-e) > \alpha(H)$ for all $e \in E(H)$. Since $H$ is triangle-free the set of vertices $N[u] := N_1(u) \cup \{u\}$ is a bipartite graph for all $u \in V(H)$. If, in addition, $H$ is an edge-critical graph then it is easily seen that $N[u]$ is a minimal destabiliser of $H$. Furthermore, we then denote $\Cr_{2}(H,N[u])$ by $\Cr_2(H,u)$. In this case there is also a unique neighbour of the apex of the stitch which is adjacent to $u$, which we name the \emph{coapex} of the stitch.

\subsubsection{Other kinds of stitches}

Suppose $M_1 = A_1 \cupdot B_1$ and $M_2 = A_2 \cupdot B_2$ are two disjoint bipartite subsets of $V(G)$. Assume moreover that if only one of $A_i$ and $B_i$ contains exactly one vertex, then it is $A_i$. We define a \emph{$C_4$ stitch-graph} $G'$ from $G$ by adding four new vertices to $G$, $\{c_1,c_2,c_3,c_4\}$, that form a cycle of length four (with $c_1c_3 \notin E(G')$). Also add edges from $c_i$ to all vertices in $A_i$ and from $c_{i+2}$ to all vertices in $B_i$ for $i \in \{1,2\}$. The resulting graph $G'$ is denoted by $\Cr_{C_4}(G;M_1,M_2)$. The quintuple $(\Cr_{C_4}(G;M_1,M_2),c_1,c_2,c_3,c_4)$ is called a \emph{$C_4$ stitch} of $G$ and $\Cr_{C_4}(G;M_1,M_2)$ is called the \emph{stitch-graph} of the $C_4$ stitch. With this definition $c_2$ is said to be the \emph{apex} of the $C_4$ stitch and $c_1$ the \emph{coapex}.

Now, instead suppose that we have three disjoint bipartite subsets of $V(G)$, say $M_1 = A_1 \cupdot B_1$, $M_2 = A_2 \cupdot B_2$ and $M_3 = A_3 \cupdot B_3$. As above we assume that if $M_i$ induces a $K_{1,\ell}$ in $G$ for some $\ell \geq 2$, then $A_i$ is the part of size one. We define a \emph{triangular stitch-graph} $G'$ from $G$ by adding the vertices $\{v_i \mid i \in [5]\}$ to $G$ and edges $v_1v_2,v_1v_3,v_1v_4,v_1v_5$ together with all edges possible between the following pairs of sets of vertices: $(\{v_2\},A_1)$, $(\{v_2\},A_2)$, $(\{v_2\},A_3)$, $(\{v_3\},B_1)$, $(\{v_4\},B_2)$ and $(\{v_5\},B_3)$. The graph $G'$ formed in such a way is denoted by $\Cr_{\bigtriangleup}(G;M_1,M_2,M_3)$. The sextuple $(\Cr_{\bigtriangleup}(G;M_1,M_2,M_3),v_1,v_2,v_3,v_4,v_5)$ is called a \emph{triangular stitch} of $G$. In this situation call $v_1$ the \emph{apex} and $v_2$ the \emph{coapex} of the triangular stitch.

For our final kind of special stitch we, as above, have three disjoint bipartite subsets, $M_1$, $M_2$ and $M_3$, subject to the same conditions. We form the \emph{$K_{2,3}$ stitch-graph} $G'$ from $G$ by adding five new vertices $\{v_i \mid i \in [5]\}$ to $G$ and edges $v_1v_3,v_1v_4,v_1v_5,v_2v_3,v_2v_4,v_2v_5$. Furthermore we add all edges between the pairs $(\{v_1\},A_1)$, $(\{v_2\},B_1)$, $(\{v_3\},A_2)$, $(\{v_4\},B_2)$, $(\{v_4\},A_3)$, $(\{v_5\},B_3)$ and $(\{v_3\},A_3)$. Denote by $\Cr_{K_{2,3}}(G; M_1,M_2,M_3)$ the graph $G'$. At last we say that the sextuple $(\Cr_{K_{2,3}}(G;M_1,M_2,M_3),v_1,v_2,v_3,v_4,v_5)$ is a \emph{$K_{2,3}$ stitch} of $G$. In this case we may call $v_1$ the \emph{apex} of the stitch. We have no need to label any of the vertices as the coapex.

In all the above notation we, as for $d$-stitches, may just write $u \in V(G)$ instead of $M_i$ in the $\Cr$-notation if $M_i = N[u]$ with $A_i = \{u\}$.

All the stitches described are part of a more general theory of stitches and patterns (studied in \cite{carc}) which we will not go into here.

\subsection{Patterned graphs} \label{patterned}

A \emph{pattern} is (or will here be taken to mean) a connected triangle-free graph of maximum valency at most three together with some directions assigned to edges that are incident to trivalent vertices, with the property that each trivalent vertex has out-valency 1 or 3 (and no edge has more than one direction). The proposition that there is a directed edge in the pattern from the vertex $p$ to the vertex $q$ we write as $p \diredge q$.

Now we go on to describe how to generate a graph (unique up to isomorphism) from a given pattern $\mathcal{P}$. It can be shown, but is not shown here (see \cite{carc} for more details), that any of the apparent choices made in the construction of the graph do not matter, since the resulting graphs will always belong to the same isomorphism class. There are also some claims of existence of certain local substructures in the recursive step of defining a graph. That these structures exist, and are unique if claimed to be, may be shown by induction, but is not done here either.

Order the vertices $V(\mathcal{P}) = \{p_1,p_2,\dots,p_n\}$ in any order.
We will recursively construct a graphs $G_\ell$ and maps
\begin{align*}
\apex_\ell &: \{p_1,\dots,p_\ell\} \to V(G_\ell) \\
\coapex_\ell &: \{p_1,\dots,p_\ell\} \to V(G_\ell) \cup \{\emptyset\} \\
\Act_\ell &: \{p_1,\dots,p_\ell\} \to \{S \subseteq V(G_\ell) \mid \abs{S} \leq 2\} \\
\faces_\ell &: \{(p_i,x) \mid p_i \in \{p_1,\dots,p_\ell\}, x \in \Act_\ell(p_i)\} \to \{p_1,\dots,p_\ell\} \cup \{\emptyset\}.
\end{align*}
The maps are used to keep track of the local structure in the graph $G_\ell$ and determines how we compute the next graph $G_{\ell + 1}$. We want to use some stitch to achieve this, and the maps determines what base sets can be used for the stitch. We call vertices in $G_\ell$ in the images of $\Act_\ell$-maps \emph{active} vertices. We say that an active $G_\ell$-vertex $x$ at some pattern vertex $p$ \emph{faces} another pattern vertex $q$ if $\faces_\ell(p,x) = q$.

Define $\mathcal{P}_k = \{p_1,\dots,p_k\}$ with $\mathcal{P}_0 = \emptyset$. Let $\deg_k(p)$ denote the valency of vertex $p$ in the induced subpattern of $\mathcal{P}$ on the vertices $\{p_1,\dots,p_k\}$.

\textbf{Case 1: }If $\deg_k(p_k) = 0$ then we recursively define the graph $G_k$ and the maps $\apex_k$, $\coapex_k$, $\Act_k$ and $\faces_k$ as follows.
\begin{enumerate}[(i)]
\item $G_k = G_{k-1} + P_2 = \Cr_1(G_{k-1}, \emptyset)$. Suppose that $\{v_1,v_2\} = V(G_k)\setminus V(G_{k-1})$.
\item $\apex_k(p_k) = v_1$, and $\apex_k|_{\mathcal{P}_{k-1}} = \apex_{k-1}$.
\item $\coapex_k(p_k) = v_2$, and $\coapex_k|_{\mathcal{P}_{k-1}} = \coapex_{k-1}$.
\item $\Act_k(p_k) = \{v_1,v_2\}$, and $\Act_k|_{\mathcal{P}_{k-1}} = \Act_{k-1}$.
\item $\faces_k(p_k,u) = \emptyset$ for both $u \in \Act_k(p_k)$.
\end{enumerate}

\textbf{Case 2: } If $\deg_k(p_k) = 1$, say $p_\ell$ is the neighbour of $p_k$ in $\mathcal{P}$ with $\ell < k$, and $\deg_k(p_\ell) = 1$ then we define the graph $G_k$ and the maps as follows.
\begin{enumerate}[(i)]
\item $G_k = \Cr_2(G_{k-1}, \coapex_{k-1}(p_\ell))$. Suppose that $v_1$ is the apex, and $v_2$ the coapex, of the 2-stitch.
\item $\apex_k(p_k) = v_1$, $\apex_k|_{\mathcal{P}_{k-1}} = \apex_{k-1}$.
\item $\coapex_k(p_k) = v_2$, and $\coapex_k|_{\mathcal{P}_{k-1}} = \coapex_{k-1}$.
\item $\Act_k(p_k) = \{v_1,v_2\}$, and $\Act_k|_{\mathcal{P}_{k-1}} = \Act_{k-1}$.
\item $\faces_k(p_k,v_1) = \emptyset$, $\faces_k(p_k,v_2) = p_\ell$, $\faces_k(p_\ell,\coapex_{k-1}(p_\ell)) = p_k$, \\
$\faces_k(p_\ell,\apex_{k-1}(p_\ell)) = \emptyset$ and $\faces_k(p,x) = \faces_{k-1}(p,x)$ for all $p \in \mathcal{P}_k \setminus \{p_k,p_\ell\}$ and $x \in \Act_{k-1}(p)$.
\end{enumerate}
Note that we then end up with a $C_5$-component in $G_k$ containing the apices and coapices of $p_k$ and $p_\ell$. The apices are at distance two in the component while the coapices are adjacent and are facing the opposite vertex in $\{p_k,p_\ell\}$ (to which it is the coapex for).

\textbf{Case 3:} If $p_k$ and $p_\ell$ are as in Case 2, except for that $\deg_k(p_\ell) > 1$. Exactly one of the following four cases then holds.
\begin{gather}
 \tag{C1} \label{C1} \left\{
 \begin{split}
 & \deg_k(p_\ell) = 2 \text{ or } \\
 & (\deg_k(p_\ell) = 3, \outdeg(p_\ell) = 1, \text{ say } p_\ell \diredge p_i,\; i < k, \text{ and } \\
 & \quad \faces_{k-1}(p_\ell, \apex_{k-1}(p_\ell)) = p_i)
 \end{split} \right. \\
 \tag{C2} \label{C2} \deg_k(p_\ell) = 3 \text{ and } \outdeg(p_\ell) = 3 \\
 \tag{C3} \label{C3} \left\{
 \begin{split}
 & \deg_k(p_\ell) = 3, \outdeg(p_\ell) = 1, \text{ say } p_\ell \diredge p_i,\; i < k, \text{ and } \\
 & \quad \faces_{k-1}(p_\ell, \coapex_{k-1}(p_\ell)) = p_i
 \end{split} \right. \\
 \tag{C4} \label{C4} \deg_k(p_\ell) = 3, \outdeg(p_\ell) = 1, p_\ell \diredge p_k
\end{gather}
We define the graph $G_k$ and the maps recursively as follows.
\begin{enumerate}[(i)]
 \item $G_k = \Cr_2(G_{k-1}, M)$, where $M$ is determined by the cases above. In case \eqref{C1} we set $M = \apex_{k-1}(p_\ell)$ (or equivalently $M = N_{G_{k-1}}[\apex_{k-1}(p_\ell)]$) and in case \eqref{C3} we take $M = \coapex_{k-1}(p_\ell)$. In cases \eqref{C2} and \eqref{C4} there is a unique 4-cycle, $C_{p_\ell}$, through $\apex_{k-1}(p_\ell)$. In case \eqref{C4} let $M = C_{p_\ell}$ be this cycle and in case \eqref{C2} let $M$ be the uniquely determined path of length four in $G_{k-1}$ that contains $\apex_{k-1}(p_\ell)$ and $\coapex_{k-1}(p_\ell)$ but no other vertex of $C_{p_\ell}$.
 \item $\apex_k(p_k) = v_1$, where $v_1$ is the apex of the 2-stitch for which $G_k$ is the stitch-graph in the definition of $G_k$. $\apex_k|_{\mathcal{P}_{k-1}} = \apex_{k-1}$.
 \item $\coapex_k(p_k) = \emptyset$ in cases \eqref{C2} and \eqref{C4} and $\coapex_k(p_k) = v_2$ where $v_2$ is the coapex of the 2-stitch, by which we defined $G_k$,  in the other two cases (which in these cases is defined). $\coapex_k|_{\mathcal{P}_{k-1}} = \coapex_{k-1}$.
 \item $\Act_k(p_k) = \{v_1\}$ in cases \eqref{C2} and \eqref{C4} but $\Act_k(p_k) = \{v_1,v_2\}$ in the other two cases. If $\deg_k(p_\ell) = 3$ then $\Act_k(p_\ell) = \emptyset$, and if furthermore $p_\ell \diredge p_i$ with $\faces_{k-1}(p_i,v) = p_\ell$ for some $v \in \Act_{k-1}(p_i)$ then $\Act_{k}(p_i) = \Act_{k-1}(p_i) \setminus \{v\}$. For all other $p \in \mathcal{P}_{k-1}$ let $\Act_k(p) = \Act_{k-1}(p)$.
 \item $\faces_k(p_k,v_1) = \emptyset$ in all four cases, and $\faces_k(p_k,v_2) = p_\ell$, in cases \eqref{C1} and \eqref{C3}.  If $\deg_k(p_\ell) \neq 2$, then $\faces_k(p,x) = \faces_{k-1}(p,x)$ for all $p \in \mathcal{P}_{k-1}$ and $x \in \Act_{k-1}(p)$. On the other hand if $\deg_k(p_\ell) = 2$, then $\faces_k(p_\ell, \apex_{k-1}(p_\ell)) = p_k$ but similarly $\faces_k(p,x) = \faces_{k-1}(p,x)$ for all $p \in \mathcal{P}_{k-1}$ and $x \in \Act_{k-1}(p)$ such that $x \neq \apex_{k-1}(p_\ell)$.
\end{enumerate}

For the remaining three cases, cases 4-6, the induced maps $f|_{\mathcal{P}_{k-1}}$ are defined recursively in the exact same way as in cases 2 and 3 (the $\deg_k(p_k) = 1$-cases), for $f \in \{\apex_k,\coapex_k,\Act_k\}$. Also $\faces_k|_{\mathcal{P}_k \times V(G_k)}$ is determined in the same way as in those cases.

\textbf{Case 4:} If $\deg_k(p_k) = 2$, say $p_{\ell_1}$ and $p_{\ell_2}$ are the neighbours of $p_k$ in $\mathcal{P}$ with $\ell_1, \ell_2 < k$. We assume that $\deg_k(p_{\ell_1}) \leq \deg_k(p_{\ell_2})$ and if $\deg_k(p_{\ell_1}) = \deg_k(p_{\ell_2}) = 3$ but only one of $p_{\ell_1}$ and $p_{\ell_2}$ has a directed edge to $p_k$, then let it be $p_{\ell_1}$, i.e. we then have $p_{\ell_1} \diredge p_k$.

We define the graph $G_k$ and the maps recursively as follows.

\begin{enumerate}[(i)]
 \item $G_k = \Cr_{C_4}(G_{k-1};M_1,M_2)$ is defined to be the stitch-graph of a $C_4$ stitch $(G_k,c_1,c_2,c_3,c_4)$ where bipartite subsets $M_1$ and $M_2$ are determined for $p_{\ell_1}$ and $p_{\ell_2}$, respectively, in the exact same way as $M$ was determined for $p_\ell$ in Case 3.
 \item $\apex_k(p_k) = c_2$.
 \item
$$\coapex_k(p_k) =
\begin{cases}
c_1, \quad & \text{if } \deg_k(p_{\ell_1}) < 3 \text{ or } \deg_k(p_{\ell_2}) < 3 \\
\emptyset, \quad & \text{otherwise.}
\end{cases}$$
 \item $\Act_k(p_k) = \emptyset$ if $\deg_k(p_{\ell_1}) = \deg_k(p_{\ell_2}) = 3$, and $p_{\ell_i} \diredge p_k$ for both $i \in [2]$. $\Act_k(p_k) = \{c_2\}$ if $\deg_k(p_{\ell_1}) = \deg_k(p_{\ell_2}) = 3$ but not both $p_{l_1}$ and $p_{l_2}$ have a directed edge to $p_k$. $\Act_k(p_k) = \{c_1\}$  if $\deg_k(p_{\ell_1}) \neq \deg_k(p_{\ell_2}) = 3$ and $p_{\ell_2} \diredge p_k$. In all other cases we set $\Act_k(p_k) = \{c_1,c_2\}$.
 \item For all $i \in [2]$ and $u \in \Act_k(p_k)$ we set
 $$\faces_k(p_k,u) = \begin{cases}
p_{\ell_i}, & \text{if } u \in N_{G_k}(M_i) \\
\emptyset, & \text{otherwise.} \end{cases}$$
\end{enumerate}

\textbf{Case 5:} If $\deg_k(p_k) = 3$ and $\outdeg(p_k) = 3$ suppose that $p_{\ell_1}$, $p_{\ell_2}$ and $p_{\ell_3}$ are the three neighbours of $p_k$ in $\mathcal{P}$. We define the graph $G_k$ and the maps by the following.
\begin{enumerate}[(i)]
 \item $G_k = \Cr_{\bigtriangleup}(G_{k-1}; x_1,x_2,x_3)$ is defined as the triangular stitch-graph to a triangular stitch $(G_k,v_1,v_2,v_3,v_4,v_5)$, where
$$x_i = \begin{cases}
\apex_{k-1}(p_{\ell_i}), & \text{if \eqref{C1}  holds for } p_{\ell_i} \\
\coapex_{k-1}(p_{\ell_i}), & \text{if \eqref{C3}  holds for } p_{\ell_i},
\end{cases}$$
for all $i \in [3]$. This is well defined, since in this case either \eqref{C1} or \eqref{C3} must hold for all the neighbours of $p_k$.
 \item $\apex_k(p_k) = v_1$, the apex of the triangular stitch.
 \item $\coapex_k(p_k) = v_2$, the coapex of the triangular stitch.
 \item $\Act_k(p_k) = \emptyset$.
 \item $\faces_k(p_k,\cdot)$ is the empty function since it has domain $\Act_k(p_k) = \emptyset$.
\end{enumerate}

\textbf{Case 6:} Finally, if $\deg_k(p_k) = 3$ and $\outdeg(p_k) = 1$ suppose that $p_{\ell_1}$, $p_{\ell_2}$ and $p_{\ell_3}$ are the three neighbours of $p_k$ in $\mathcal{P}$, and that $p_k \diredge p_{\ell_3}$. We then define the graph $G_k$ and the maps by the following.
\begin{enumerate}
 \item $G_k = \Cr_{K_{2,3}}(G_{k-1}; M_1,M_2,x_3)$ is defined to be the stitch-graph of a $K_{2,3}$ stitch $(G_k,v_1,v_2,v_3,v_4,v_5)$, where $M_1$ and $M_2$ are bipartite subsets that are determined for $p_{\ell_1}$ and $p_{\ell_2}$, respectively, in the same way as $M$ was determined for $p_\ell$ in Case 3. The vertex $x_3$ is determined in the same way as in Case 5.
 \item $\apex_k(p_k) = v_1$, the apex of the $K_{2,3}$ stitch.
 \item $\coapex_k(p_k) = \emptyset$.
 \item $\Act_k(p_k) = \emptyset$.
 \item $\faces_k(p_k,\cdot)$ is the empty function since it has domain $\Act_k(p_k) = \emptyset$. 
\end{enumerate}

This completes the recursive definition of $G_k$ and the maps $\apex_k$, $\coapex_k$, $\Act_k$ and $\faces_k$.

We will now remark on how the number of vertices, number of edges and the independence number of $G  = G_n$ is determined from the pattern $\mathcal{P}$ on $\{p_1,\dots,p_n\}$. By showing that each of the stitches in the recursive construction of $G$ are based at minimal destabilisers it can be proven that $\alpha(G_k) = \alpha(G_{k-1}) + 1$. Hence $\alpha(G) = n = n(\mathcal{P})$, the number of vertices in the pattern. Note also that none of the stitches may introduce triangles into the graph and therefore $G$ is triangle-free.

It is not difficult to verify that $n(G) = 2n + e(\mathcal{P})$ and $e(G) = n + 2e(\mathcal{P}) + \frac{1}{2} \sum_{p \in V(\mathcal{P})} \deg(p)^2$. The first identity can be verified by checking that $n(G_k) - n(G_{k-1}) = 2 + \deg_k(p_k)$ in all of the six cases. The second identity can be shown by first showing that in cases \eqref{C1} and \eqref{C3} with $\deg_k(p_\ell) = 3$ we must have that $\apex_{k-1}(p_\ell)$ and $\coapex_{k-1}(p_\ell)$, respectively, are trivalent in $G_{k-1}$. Then one verifies the identity by showing that $e(G_k) - e(G_{k-1}) = 1 + \sum_{p \in N(p_k)} \deg_k(p) + \frac{1}{2} \deg_k(p_k)^2 + \frac{3}{2} \deg_k(p_k)$ in each of the six cases.

\subsubsection{$H_{13}$-decorated patterned graphs}

A $H_{13}$-decorated pattern (or just $H_{13}$-pattern) $(\mathcal{P},\mathcal{H})$ consists of an ordinary pattern $\mathcal{P}$ and a set of hyperedges $\mathcal{H} \subseteq \binom{V(\mathcal{P})}{4}$ of size four such that:
\begin{enumerate}[(i)]
 \item If $H \in \mathcal{H}$ then $H$ induces a $C_4$ in $\mathcal{P}$.
 \item If $H \in \mathcal{H}$ then $e(H,V(\mathcal{P})\setminus H) \leq 1$.
 \item If $H \in \mathcal{H}$ contains a vertex $v$ which is trivalent in $\mathcal{P}$, then $\outdeg(v) = 1$.
\end{enumerate}
The second condition can be relaxed to get more general patterns. These will however not be considered here. These patterns are named after the unique $(3,5; 13,26)$-graph, which is of special interest withing Ramsey theory and is commonly known as $H_{13}$. $H_{13}$ can be described as the circulant graph on vertices $\{v_0,v_1,\dots,v_{12}\}$ such that all there is an edge between $v_i$ and $v_j$ if their indices differ by a number which is congruent to either 1 or 5 modulo 13. This graph is illustrated in Figure \ref{H13}.

We will now give a description of how to obtain a $H_{13}$-patterned graph from a $H_{13}$-pattern. First form an \emph{auxiliary ordinary pattern} $\mathcal{P}'$ from $(\mathcal{P},\mathcal{H})$ by contracting each $H \in \mathcal{H}$ to a single vertex $v_H$. If there is an edge directed outward from $H$ in $\mathcal{P}$, then we make this edge undirected in $\mathcal{P}'$. This forms an ordinary pattern on $k = n(\mathcal{P}) - 3|\mathcal{H}|$ vertices. We may then recursively compute the patterned graph $G_k$ from pattern $\mathcal{P}'$ by the procedure described above.

Let $u_{H,1}$ and $u_{H,2}$ be the apex and coapex of $v_H$ in $G_k$, respectively. If the coapex is not defined (i.e. $\coapex_k(v_H) = \emptyset$) take $u_{H,2}$ to be any neighbour of $u_{H,1}$ in $G_k$. The idea is now to make a $H_{13}$-decoration at each of the vertices $v_H$ by adding for each such vertex 11 new vertices $w_{H,1},\dots,w_{H,11}$ to $G_k$, forming the new graph $G$, in such a way that
\begin{enumerate}[(i)]
 \item $\{u_{H,1},u_{H,2},w_{H,3},w_{H,4},\dots,w_{H,13}\}$ induce a $H_{13}$ in $G$, in the manner as specified in Figure \ref{H13}.
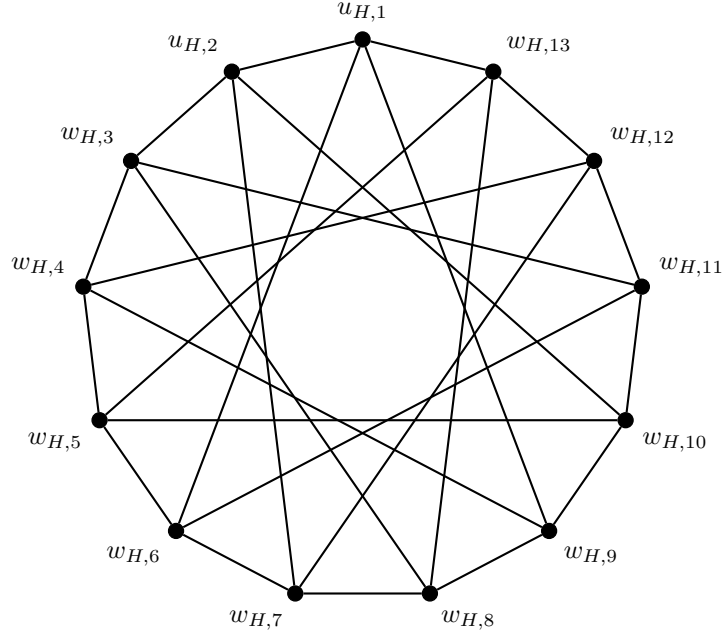
\begin{figure}
\begin{center}
\begin{tikzpicture}[scale=1.7,rotate=90]
  \GraphInit[vstyle=Classic]
  \renewcommand*{\VertexSmallMinSize}{4pt}
  \Vertices[Math,unit=2.2,Lpos=90]{circle}{u_{H,1},u_{H,2},w_{H,3},w_{H,4},w_{H,5},w_{H,6},w_{H,7},w_{H,8},w_{H,9},w_{H,10},w_{H,11},w_{H,12},w_{H,13}}
  \Edge(u_{H,1})(u_{H,2})
  \Edge(u_{H,1})(w_{H,6})
  \Edge(u_{H,1})(w_{H,13})
  \Edge(u_{H,1})(w_{H,9})
  \Edge(u_{H,2})(w_{H,3})
  \Edge(u_{H,2})(w_{H,10})
  \Edge(u_{H,2})(w_{H,7})
  \Edge(w_{H,6})(w_{H,11})
  \Edge(w_{H,6})(w_{H,7})
  \Edge(w_{H,6})(w_{H,5})
  \Edge(w_{H,13})(w_{H,12})
  \Edge(w_{H,13})(w_{H,8})
  \Edge(w_{H,13})(w_{H,5})
  \Edge(w_{H,9})(w_{H,4})
  \Edge(w_{H,9})(w_{H,10})
  \Edge(w_{H,9})(w_{H,8})
  \Edge(w_{H,3})(w_{H,4})
  \Edge(w_{H,3})(w_{H,11})
  \Edge(w_{H,3})(w_{H,8})
  \Edge(w_{H,4})(w_{H,12})
  \Edge(w_{H,4})(w_{H,5})
  \Edge(w_{H,11})(w_{H,12})
  \Edge(w_{H,11})(w_{H,10})
  \Edge(w_{H,12})(w_{H,7})
  \Edge(w_{H,10})(w_{H,5})
  \Edge(w_{H,7})(w_{H,8})
\end{tikzpicture}
\end{center}
\caption{The induced graph $H_{13}$, which is the unique $(3,5;13,26)$-graph.} \label{H13}
\end{figure}

 \item If there is an directed edge outward from $H$ in $\mathcal{P}$: For every vertex $w \in N_{G_k}(u_{H,1})$ there is an edge $w_{H,7}w \in E(G)$, and for every vertex $w \in N_{G_k}(u_{H,2})$ there are edges $w_{H,6}w \in E(G)$ and $w_{H,9}w \in E(G)$. \label{tva}
 \item If there is no directed edge outward from $H$ in $\mathcal{P}$: For every vertex $w \in N_{G_k}(u_{H,2})$ there are edges $w_{H,6}w,w_{H,13}w,w_{H,9}w \in E(G)$. \label{tre}
 \item $e(G)$ is as small as possible satisfying the above conditions.
\end{enumerate}

The case (\ref{tva}) from above corresponds to a stitch from the $H_{13}$ based at a bipartite minimal destabiliser which induces a $K_{2,3}^-$-graph ($K_{2,3}^-$ is the graph obtained from $K_{2,3}$ by removing one edge). The other possibility, (\ref{tre}), corresponds to a stitch from the $H_{13}$ based at $N[u_{H,1}]$.

One can verify that now we get $n(G) = 2n(\mathcal{P}) + e(\mathcal{P}) + |\mathcal{H}|$, $\alpha(G) = n(\mathcal{P})$. Moreover, one can show that each $H_{13}$-patterned graph corresponds uniquely up to isomorphism to a $H_{13}$-pattern. The description of $H_{13}$-patterned graphs as described here may seem quite artificial. One can fit these more naturally into a general framework of patterned graphs (see \cite{carc}), but this complicates the description somewhat and is therefore not done in this report.

\section{Computation}

To obtain the values in Table \ref{maintable} and the tables in the appendix we used a computer program written in C which heavily relies on the \texttt{nauty} package (see \cite{webnauty} and \cite{articlenauty}) by B. McKay.

The program works by first by using the \texttt{geng} utility from \texttt{nauty} to first generate all triangle-free connected graphs with maximum valency at most three. Then we generate all auxiliary patterns with these graphs as their underlying graphs (the graphs we get by removing directions from edges in a pattern). Using \texttt{nauty}:s canonical labelling one can make sure that one gets a list of only non-isomorphic patterns.

From each $H_{13}$-pattern obtained in such a way we compute the unique corresponding $H_{13}$-patterned graph by following the recursive procedure from Section \ref{patterned}.

The result of this computation of all patterned $(3,k)$-graphs and the C-code of the program that was used to compute them are available upon request to the author.

We, in Appendix \ref{apxlist}, count the number of $H_{13}$-patterned graphs for $(3,k)$, $k \leq 12$ obtained in this way.

\section{Example: Patterned and non-patterned $(3,6;16,32)$- and $(3,6;16,33)$-graphs}

Since the general abstract description above is quite technical we here give example calculations of how one obtains the four $H_{13}$-patterned $(3,6;16)$-graphs. We start by determining the possible underlying graphs of patterns that correspond to $(3,6;16)$-graphs.

If $(\mathcal{P},\mathcal{H})$ is an $H_{13}$-pattern that corresponds to a $(3,6;16)$-graph, $G$, then $k = 5 = n(\mathcal{P})$. So in particular $\mathcal{H}$ contains 0 or 1 hyperedge (since the pattern can not contain more than one cycle of length four). The underlying graph of the pattern (which in this case is an ordinary pattern) should contain $\frac{1}{2}(n(G) - n(\mathcal{P})) = 6$ edges if $|\mathcal{H}| = 0$ and 5 edges if $|\mathcal{H}| = 1$. There is only one connected triangle-free graph with maximum valency at most three with five vertices and four edges, namely $K_{2,3}$. This is therefore the only possible underlying graph of an ordinary pattern, i.e. if $|\mathcal{H}| = 0$. In the case that $|\mathcal{H}| = 1$ the underlying graph of the pattern should have five edges and a cycle of length four, the only possible graph is then $K_{2,3}^-$.

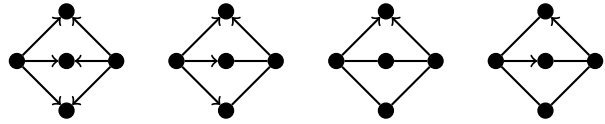
\begin{figure}
\begin{center}
\begin{tikzpicture}[scale=0.3]
  \GraphInit[vstyle=Classic]
  \renewcommand*{\VertexSmallMinSize}{2pt}
  \Vertices[Math,unit=2.2,NoLabel]{circle}{v1,u1,v2,u3}
  \Vertex[Math,NoLabel,x=0,y=0]{u2}
  \Edge[style=<-](u1)(v1)
  \Edge[style=<-](u1)(v2)
  \Edge[style=<-](u2)(v1)
  \Edge[style=<-](u2)(v2)
  \Edge[style=<-](u3)(v1)
  \Edge[style=<-](u3)(v2)
\end{tikzpicture} $\quad$
\begin{tikzpicture}[scale=0.3]
  \GraphInit[vstyle=Classic]
  \renewcommand*{\VertexSmallMinSize}{2pt}
  \Vertices[Math,unit=2.2,NoLabel]{circle}{v1,u1,v2,u3}
  \Vertex[Math,NoLabel,x=0,y=0]{u2}
  \Edge[style=<-](u1)(v1)
  \Edge[style=<-](u1)(v2)
  \Edge(u2)(v1)
  \Edge[style=<-](u2)(v2)
  \Edge(u3)(v1)
  \Edge[style=<-](u3)(v2)
\end{tikzpicture} $\quad$
\begin{tikzpicture}[scale=0.3]
  \GraphInit[vstyle=Classic]
  \renewcommand*{\VertexSmallMinSize}{2pt}
  \Vertices[Math,unit=2.2,NoLabel]{circle}{v1,u1,v2,u3}
  \Vertex[Math,NoLabel,x=0,y=0]{u2}
  \Edge[style=<-](u1)(v1)
  \Edge[style=<-](u1)(v2)
  \Edge(u2)(v1)
  \Edge(u2)(v2)
  \Edge(u3)(v1)
  \Edge(u3)(v2)
\end{tikzpicture} $\quad$
\begin{tikzpicture}[scale=0.3]
  \GraphInit[vstyle=Classic]
  \renewcommand*{\VertexSmallMinSize}{2pt}
  \Vertices[Math,unit=2.2,NoLabel]{circle}{v1,u1,v2,u3}
  \Vertex[Math,NoLabel,x=0,y=0]{u2}
  \Edge[style=<-](u1)(v1)
  \Edge(u1)(v2)
  \Edge(u2)(v1)
  \Edge[style=<-](u2)(v2)
  \Edge(u3)(v1)
  \Edge(u3)(v2)
\end{tikzpicture}
\end{center}
\caption{The four possible patterns where $|\mathcal{H}| = 0$.} \label{patternsH0}
\end{figure}

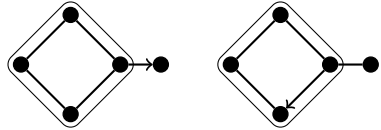
\begin{figure}
\begin{center}
\begin{tikzpicture}[scale=0.3]
  \GraphInit[vstyle=Classic]
  \renewcommand*{\VertexSmallMinSize}{2pt}
  \Vertices[Math,unit=2.2,NoLabel]{circle}{v1,u1,v2,u3}
  \Vertex[Math,NoLabel,x=4,y=0]{u2}
  \Edge(u1)(v1)
  \Edge(u1)(v2)
  \Edge[style=<-](u2)(v1)
  \Edge(u3)(v1)
  \Edge(u3)(v2)
  \draw[rounded corners,rotate=45] (-2.1, -2.1) rectangle (2.1, 2.1) {};
\end{tikzpicture} $\quad$
\begin{tikzpicture}[scale=0.3]
  \GraphInit[vstyle=Classic]
  \renewcommand*{\VertexSmallMinSize}{2pt}
  \Vertices[Math,unit=2.2,NoLabel]{circle}{v1,u1,v2,u3}
  \Vertex[Math,NoLabel,x=4,y=0]{u2}
  \Edge(u1)(v1)
  \Edge(u1)(v2)
  \Edge(u2)(v1)
  \Edge[style=<-](u3)(v1)
  \Edge(u3)(v2)
  \draw[rounded corners,rotate=45] (-2.1, -2.1) rectangle (2.1, 2.1) {};
\end{tikzpicture}
\end{center}
\caption{The two possible patterns where $|\mathcal{H}| = 1$.} \label{patternsH1}
\end{figure}

In the case $|\mathcal{H}| = 0$ we have four distinct ways to make a pattern $\mathcal{P}$ with underlying graph $K_{2,3}$. These are illustrated in Figure \ref{patternsH0}. If $|\mathcal{H}| = 1$ we instead have two $H_{13}$-patterns, which are illustrated in Figure \ref{patternsH1}.

Let us now further specialise this example calculation to a pair of patterns, one for $|\mathcal{H}| = 0$ and one for $|\mathcal{H}| = 1$ and illustrate in steps how one obtains graphs from these patterns.

\begin{figure}
\begin{center}
\begin{tikzpicture}[scale=0.3]
  \GraphInit[vstyle=Classic]
  \renewcommand*{\VertexSmallMinSize}{2pt}
  \Vertices[Math,unit=2.2]{circle}{p_5,p_2,p_1,p_4}
  \Vertex[Math,x=0,y=0,Lpos=-90]{p_3}
  \Edge[style=<-](p_2)(p_5)
  \Edge[style=<-](p_2)(p_1)
  \Edge(p_3)(p_5)
  \Edge[style=<-](p_3)(p_1)
  \Edge(p_4)(p_5)
  \Edge[style=<-](p_4)(p_1)
\end{tikzpicture}
\end{center}
\caption{An example of a pattern with the vertices ordered.} \label{orderedH0}
\end{figure}
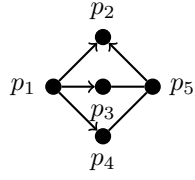

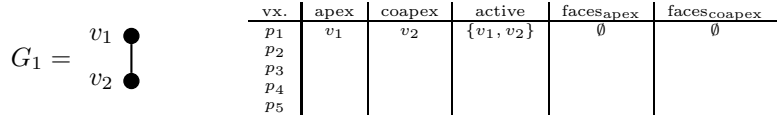
\begin{figure}[H]
\begin{center}
\begin{minipage}{0.25\textwidth}
\begin{tikzpicture}[scale=0.3]
  \GraphInit[vstyle=Classic]
  \renewcommand*{\VertexSmallMinSize}{2pt}
  \Vertex[Math,x=0,y=1,Lpos=-180]{v_1}
  \Vertex[Math,x=0,y=-1,Lpos=-180]{v_2}
  \Edge(v_1)(v_2)
  \node at (-4,0) {$G_1 = $};
\end{tikzpicture}
\end{minipage}
\begin{minipage}{0.6\textwidth}
\tiny{
\begin{tabular}{c | c | c | c | c | c}
  vx. & $\apex$ & $\coapex$ & $\Act$ & $\faces_{\apex}$ & $\faces_{\coapex}$ \\
  \hline
  $p_1$ & $v_1$ & $v_2$ & $\{v_1,v_2\}$ & $\emptyset$ & $\emptyset$ \\
  $p_2$ & & & & & \\
  $p_3$ & & & & & \\
  $p_4$ & & & & & \\
  $p_5$ & & & & &
\end{tabular}
}
\end{minipage}
\end{center}
\caption{k = 1, we use Case 1 in the recursive construction since $\deg_1(p_1) = 0$.} \label{step1}
\end{figure}

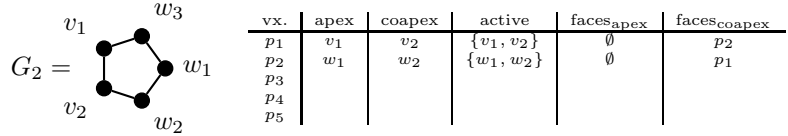
\begin{figure}[H]
\begin{center}
\begin{minipage}{0.25\textwidth}
\begin{tikzpicture}[scale=0.3]
  \GraphInit[vstyle=Classic]
  \renewcommand*{\VertexSmallMinSize}{2pt}
  \Vertices[Math,unit=1.5]{circle}{w_1,w_3,v_1,v_2,w_2}
  \Edges(w_1,w_3,v_1,v_2,w_2,w_1)
  \node at (-4,0) {$G_2 = $};
\end{tikzpicture}
\end{minipage}
\begin{minipage}{0.6\textwidth}
\tiny{
\begin{tabular}{c | c | c | c | c | c}
  vx. & $\apex$ & $\coapex$ & $\Act$ & $\faces_{\apex}$ & $\faces_{\coapex}$ \\
  \hline
  $p_1$ & $v_1$ & $v_2$ & $\{v_1,v_2\}$ & $\emptyset$ & $p_2$ \\
  $p_2$ & $w_1$ & $w_2$ & $\{w_1,w_2\}$ & $\emptyset$ & $p_1$ \\
  $p_3$ & & & & & \\
  $p_4$ & & & & & \\
  $p_5$ & & & & &
\end{tabular}
}
\end{minipage}
\end{center}
\caption{k = 2, we use Case 2 in the recursive construction since $\deg_2(p_2) = \deg_2(p_1) = 1$.} \label{step2}
\end{figure}

In the case where $|\mathcal{H}| = 0$ we take the pattern which contains exactly four directed edges (the second one from the left in Figure \ref{patternsH0}) as an example. We choose some order $p_1,p_2,p_3,p_4,p_5$ of the vertices in the pattern; as in Figure \ref{orderedH0}. We start by having $\mathcal{P}_0 = \emptyset$ and $G_0 = \emptyset$. We now step through the recursive construction for $k = 1,2,3,4$ and $5$ in this example. Each step with graphs $G_k$ and corresponding maps are illustrated in Figures \ref{step1} through \ref{step5}.

\begin{figure}[H]
\begin{center}
\begin{minipage}{0.25\textwidth}
\begin{tikzpicture}[scale=0.5]
  \GraphInit[vstyle=Classic]
  \renewcommand*{\VertexSmallMinSize}{2pt}
  \Vertex[Math,x=0,y=0,Lpos=-180]{u_1}
  \Vertex[Math,x=0,y=1,Lpos=-180]{u_2}
  \Vertex[Math,x=1,y=0,Lpos=-90]{u_3}
  \Vertex[Math,x=1,y=1,Lpos=90]{v_1}
  \Vertex[Math,x=2,y=0,Lpos=-90]{v_2}
  \Vertex[Math,x=2,y=1,Lpos=90]{w_3}
  \Vertex[Math,x=3,y=0,Lpos=-90]{w_2}
  \Vertex[Math,x=3,y=1,Lpos=90]{w_1}
  \Edges(u_1,u_2,v_1,w_3,w_1,w_2,v_2,u_3,u_1)
  \Edge(v_1)(v_2)
  \Edge(u_3)(w_3)
  \node at (-2,0.5) {$G_3 = $};
\end{tikzpicture}
\end{minipage}
\begin{minipage}{0.6\textwidth}
\tiny{
\begin{tabular}{c | c | c | c | c | c}
  vx. & $\apex$ & $\coapex$ & $\Act$ & $\faces_{\apex}$ & $\faces_{\coapex}$ \\
  \hline
  $p_1$ & $v_1$ & $v_2$ & $\{v_1,v_2\}$ & $p_3$ & $p_2$ \\
  $p_2$ & $w_1$ & $w_2$ & $\{w_1,w_2\}$ & $\emptyset$ & $p_1$ \\
  $p_3$ & $u_1$ & $u_2$ & $\{u_1,u_2\}$& $\emptyset$ & $p_1$ \\
  $p_4$ & & & & & \\
  $p_5$ & & & & &
\end{tabular}
}
\end{minipage}
\end{center}
\caption{k = 3, we use Case 3 (\ref{C1}) in the recursive construction since $\deg_3(p_3) = 1 < \deg_3(p_1) = 2$.} \label{step3}
\end{figure}
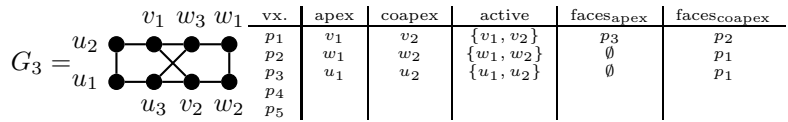

\begin{figure}[H]
\begin{center}
\begin{minipage}{0.3\textwidth}
\begin{tikzpicture}[scale=0.5]
  \GraphInit[vstyle=Classic]
  \renewcommand*{\VertexSmallMinSize}{2pt}
  \Vertex[Math,x=0,y=0,Lpos=-180]{u_1}
  \Vertex[Math,x=0,y=1,Lpos=-180]{u_2}
  \Vertex[Math,x=1,y=0,Lpos=-90]{u_3}
  \Vertex[Math,x=1,y=1,Lpos=90]{v_1}
  \Vertex[Math,x=2,y=0,Lpos=-45,Ldist=-5]{v_2}
  \Vertex[Math,x=2,y=1,Lpos=90]{w_3}
  \Vertex[Math,x=3,y=0,Lpos=-45,Ldist=-5]{w_2}
  \Vertex[Math,x=3,y=1,Lpos=90]{w_1}
  \Vertex[Math,x=0.75,y=-1.5,Lpos=-90]{x_2}
  \Vertex[Math,x=2.25,y=-1.5,Lpos=-90]{x_3}
  \Vertex[Math,x=1.5,y=-2,Lpos=-90]{x_1}
  \Edges(x_2,x_1,x_3)
  \Edges(u_2,x_2,v_2)
  \Edges(v_1,x_3,w_2)
  \Edges(u_1,u_2,v_1,w_3,w_1,w_2,v_2,u_3,u_1)
  \Edge(v_1)(v_2)
  \Edge(u_3)(w_3)
  \node at (-2,0.5) {$G_4 = $};
\end{tikzpicture}
\end{minipage}
\begin{minipage}{0.6\textwidth}
\tiny{
\begin{tabular}{c | c | c | c | c | c}
  vx. & $\apex$ & $\coapex$ & $\Act$ & $\faces_{\apex}$ & $\faces_{\coapex}$ \\
  \hline
  $p_1$ & $v_1$ & $v_2$ & $\emptyset$ & $p_1$ & $p_2$ \\
  $p_2$ & $w_1$ & $w_2$ & $\{w_1\}$ & $p_3$ & $p_1$ \\
  $p_3$ & $u_1$ & $u_2$ & $\{u_1\}$& $\emptyset$ & $p_1$ \\
  $p_4$ & $x_1$ & $\emptyset$ & $\{x_1\}$ & $\emptyset$ & \\
  $p_5$ & & & & &
\end{tabular}
}
\end{minipage}
\end{center}
\caption{k = 4, we use Case 3 (\ref{C2}) in the recursive construction since $\deg_4(p_4) = 1 < \deg_4(p_4) = 3$ and $\outdeg(p_1) = 1$.} \label{step4}
\end{figure}
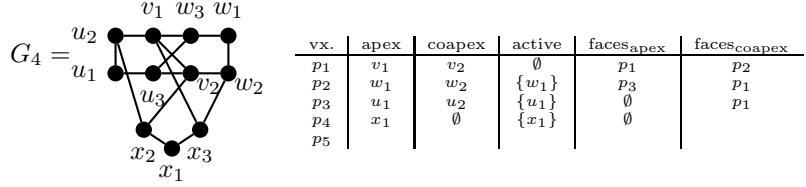

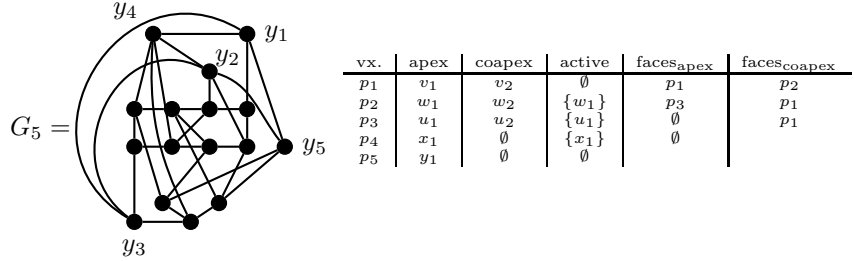
\begin{figure}[H]
\begin{center}
\begin{minipage}{0.35\textwidth}
\begin{tikzpicture}[scale=0.5]
  \GraphInit[vstyle=Classic]
  \renewcommand*{\VertexSmallMinSize}{2pt}
  \Vertex[Math,x=0,y=0,NoLabel]{u_1}
  \Vertex[Math,x=0,y=1,NoLabel]{u_2}
  \Vertex[Math,x=1,y=0,NoLabel]{u_3}
  \Vertex[Math,x=1,y=1,NoLabel]{v_1}
  \Vertex[Math,x=2,y=0,NoLabel]{v_2}
  \Vertex[Math,x=2,y=1,NoLabel]{w_3}
  \Vertex[Math,x=3,y=0,NoLabel]{w_2}
  \Vertex[Math,x=3,y=1,NoLabel]{w_1}
  \Vertex[Math,x=0.75,y=-1.5,NoLabel]{x_2}
  \Vertex[Math,x=2.25,y=-1.5,NoLabel]{x_3}
  \Vertex[Math,x=1.5,y=-2,NoLabel]{x_1}
  \Vertex[Math,x=2,y=2,Lpos=45,Ldist=-5]{y_2}
  \Vertex[Math,x=3,y=3]{y_1}
  \Vertex[Math,x=0,y=-2,Lpos=-90]{y_3}
  \Vertex[Math,x=4,y=0]{y_5}
  \Vertex[Math,x=0.5,y=3,Lpos=135]{y_4}
  \Edges(u_2,y_4,u_3)
  \Edge[style={bend right,in=195,out=-15}](y_4)(x_1)
  \Edges(x_2,y_5,x_3)
  \Edges(u_1,y_3,x_1)
  \Edges(w_3,y_2,w_2)
  \Edge(y_1)(w_1)
  \Edges(x_2,x_1,x_3)
  \Edges(u_2,x_2,v_2)
  \Edges(v_1,x_3,w_2)
  \Edges(u_1,u_2,v_1,w_3,w_1,w_2,v_2,u_3,u_1)
  \Edge(v_1)(v_2)
  \Edge(u_3)(w_3)
  \Edge(y_5)(y_1)
  \Edge[style={bend right,out=-5}](y_5)(y_2)
  \Edge(y_4)(y_1)
  \Edge(y_4)(y_2)
  \Edge[style={bend right,in=90,out=90,distance=3.7cm}](y_3)(y_1)
  \Edge[style={bend right,in=90,out=80,distance=2.5cm}](y_3)(y_2)
  \node at (-2.5,0.5) {$G_5 = $};
\end{tikzpicture}
\end{minipage}
\begin{minipage}{0.55\textwidth}
\tiny{
\begin{tabular}{c | c | c | c | c | c}
  vx. & $\apex$ & $\coapex$ & $\Act$ & $\faces_{\apex}$ & $\faces_{\coapex}$ \\
  \hline
  $p_1$ & $v_1$ & $v_2$ & $\emptyset$ & $p_1$ & $p_2$ \\
  $p_2$ & $w_1$ & $w_2$ & $\{w_1\}$ & $p_3$ & $p_1$ \\
  $p_3$ & $u_1$ & $u_2$ & $\{u_1\}$& $\emptyset$ & $p_1$ \\
  $p_4$ & $x_1$ & $\emptyset$ & $\{x_1\}$ & $\emptyset$ & \\
  $p_5$ & $y_1$ & $\emptyset$ & $\emptyset$ & &
\end{tabular}
}
\end{minipage}
\end{center}
\caption{k = 5, we use Case 6 in the recursive construction since $\deg_5(p_5) = 3$ and $\outdeg(p_5) = 1$.} \label{step5}
\end{figure}

The graph $G_5$ in Figure \ref{step5} is then the $(3,6;16,32)$-graph with the pattern in Figure \ref{orderedH0}.

\begin{figure}[H]
\begin{center}
\begin{tikzpicture}[scale=0.6]
  \GraphInit[vstyle=Classic]
  \renewcommand*{\VertexSmallMinSize}{4pt} \Vertices[Math,unit=2.2,NoLabel]{circle}{v_1,v_2,w_{H,3},w_{H,4},w_{H,5},w_{H,6},w_{H,7},w_{H,8},w_{H,9},w_{H,10},w_{H,11},w_{H,12},w_{H,13}}
  \Vertex[Math,x=-3.5,y=0,Lpos=180]{w_1}
  \Vertex[Math,x=-3,y=0.75,Lpos=135,Ldist=-4]{w_2}
  \Vertex[Math,x=-3,y=-0.75,Lpos=-135,Ldist=-4]{w_3}
  \Edges(w_{H,7},w_2,w_1,w_3,w_{H,8})
  \Edge(v_1)(v_2)
  \Edge(v_1)(w_{H,6})
  \Edge(v_1)(w_{H,13})
  \Edge(v_1)(w_{H,9})
  \Edge(v_2)(w_{H,3})
  \Edge(v_2)(w_{H,10})
  \Edge(v_2)(w_{H,7})
  \Edge(w_{H,6})(w_{H,11})
  \Edge(w_{H,6})(w_{H,7})
  \Edge(w_{H,6})(w_{H,5})
  \Edge(w_{H,13})(w_{H,12})
  \Edge(w_{H,13})(w_{H,8})
  \Edge(w_{H,13})(w_{H,5})
  \Edge(w_{H,9})(w_{H,4})
  \Edge(w_{H,9})(w_{H,10})
  \Edge(w_{H,9})(w_{H,8})
  \Edge(w_{H,3})(w_{H,4})
  \Edge(w_{H,3})(w_{H,11})
  \Edge(w_{H,3})(w_{H,8})
  \Edge(w_{H,4})(w_{H,12})
  \Edge(w_{H,4})(w_{H,5})
  \Edge(w_{H,11})(w_{H,12})
  \Edge(w_{H,11})(w_{H,10})
  \Edge(w_{H,12})(w_{H,7})
  \Edge(w_{H,10})(w_{H,5})
  \Edge(w_{H,7})(w_{H,8})
  \Edge[style={bend right,out=-80,in=-90,distance=3cm}](w_3)(v_1)
  \Edge[style={bend right}](w_2)(w_{H,13})
  \Edge[style={bend left,out=80,in=110}](w_2)(w_{H,3})
  \node at (-2.3,0.9) {$v_2$};
  \node at (-2.3,-0.9) {$v_1$};
\end{tikzpicture}
\end{center}
\caption{A $(3,6;16,33)$-graph obtained from a $H_{13}$-decorated pattern.} \label{H13stylie}
\end{figure}
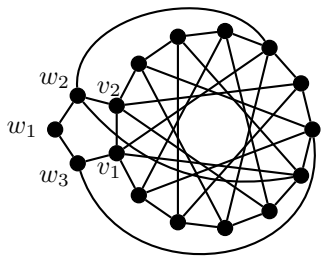

Now let us consider one of the examples where $|\mathcal{H}| = 1$ as well. Take the leftmost pattern in Figure \ref{patternsH1}.  Its auxiliary ordinary pattern $\mathcal{P}'$ is then just $K_2$. By the same procedure as before we then get that $G_2$, and the corresponding maps, is as illustrated in Figure 6. Suppose that $p_1$ was the edge formed by the contraction of the hyperedge. We then get the graph of the $H_{13}$-decorated pattern by adding eleven vertices and edges incident to them as in Figure \ref{H13stylie}. Since there is an edge out from the hyperedge in the $H_{13}$-pattern we use part (ii) in the construction to decide what edges to add between the vertices $w_i$ and the vertices that belong to the induced $H_{13}$-part of the graph.

%\printbibliography

\newpage

\appendix
\section{Tables of counts of $H_{13}$-patterned graphs}
\label{apxlist}

\thispagestyle{plain}
\markright{A: TABLES OF COUNTS OF $H_{13}$-PATTERNED GRAPHS}
\markleft{A: TABLES OF COUNTS OF $H_{13}$-PATTERNED GRAPHS}

Here follows a list of tables containing counts of the number of $H_{13}$-patterned $(3,k)$-graph found by the computer program for $k \leq 12$. A number is written in \textbf{bold} if it is the same number as can be found in the lists of McKay (at \cite{McKayMinimalRamsey}) or in the lists by Goedgebeur and Radziszowski (in \cite{HouseofGraphsMinimalRamsey} or \cite{goed-radz}). These indicate where all the graphs are $H_{13}$-patterned.
Moreover, Radziszowski and Kreher characterised all the minimal $(3,k;n)$-graphs (first for $n \leq 5(k-1)/2$ in
and for $5(k-1)/2 < n \leq 3(k-1)$ in \cite{radziszowski-kreher91}) and they are all patterned. Furthermore, if $n > 3(k-1)$ then $e(3,k;n) \geq 6n - 13k$, and the graphs achieving this with equality are all patterned (as shown in \cite{jart}). Also minimal graphs that have been characterised by any of these results have been highlighted in the tables.

For $k \leq 7$, and selected other values, we also write a subscript $(m)$ for a number in the table. Here $m$ is the total number of $(3,k)$-graphs ($H_{13}$-patterned or not) with the indicated number of vertices and edges (taken from the lists mentioned above). The subscript is omitted if all graphs are $H_{13}$-patterned. For the counts of the total number of graphs for $k \leq 8$ we refer the reader to \cite{goed-radz}.

\begin{center}
\begin{table}[h!]
\tiny
\begin{tabular}{c| c c}
edges & \multicolumn{2}{c}{number of vertices $n$} \\
$e$ & 4 & 5 \\
\hline
2 & \textbf{1} & \\
5 & & \textbf{1}
\end{tabular}
\caption{Number of $H_{13}$-patterned $(3,3)$-graphs.}
\end{table}
\end{center}

\begin{center}
\begin{table}
\tiny
\begin{tabular}{c| c c c}
edges & \multicolumn{3}{c}{number of vertices $n$} \\
$e$ & 6 & 7 & 8 \\
\hline
3 & \textbf{1} & & \\
6 & & \textbf{1} & \\
10 & & & \textbf{1}
\end{tabular}
\caption{Number of $H_{13}$-patterned $(3,4)$-graphs.}
\end{table}
\end{center}

\begin{center}
\begin{table}
\tiny
\begin{tabular}{c| c c c c c c}
edges & \multicolumn{6}{c}{number of vertices $n$} \\
$e$ & 8 & 9 & 10 & 11 & 12 & 13 \\
\hline
4 & \textbf{1} & & & & & \\
7 & & \textbf{1} & & & & \\
10 & & & \textbf{1} & & & \\
11 & & & $1_{(2)}$ & & & \\
15 & & & & \textbf{1} & & \\
16 & & & & $2_{(6)}$ & & \\
20 & & & & & \textbf{1} & \\
26 & & & & & & \textbf{1}
\end{tabular}
\caption{Number of $H_{13}$-patterned $(3,5)$-graphs.}
\end{table}
\end{center}

\begin{center}
\begin{table}
\tiny
\begin{tabular}{c| c c c c c c c}
edges & \multicolumn{7}{c}{number of vertices $n$} \\
$e$ & 10 & 11 & 12 & 13 & 14 & 15 & 16 \\
\hline
5 & \textbf{1} & & & & & & \\
8 & & \textbf{1} & & & & & \\
11 & & & \textbf{1} & & & & \\
12 & & & $1_{(4)}$ & & & & \\
15 & & & & \textbf{1} & & & \\
16 & & & & $1_{(4)}$ & & & \\
17 & & & & $2_{(45)}$ & & & \\
20 & & & & & \textbf{1} & & \\
21 & & & & & $4_{(16)}$ & & \\
25 & & & & & & \textbf{1} & \\
26 & & & & & & $3_{(7)}$ & \\
27 & & & & & & $1_{(101)}$ & \\
32 & & & & & & & $4_{(5)}$ \\
33 & & & & & & & $2_{(39)}$
\end{tabular}
\caption{Number of $H_{13}$-patterned $(3,6)$-graphs.}
\end{table}
\end{center}

\begin{center}
\begin{table}
\tiny
\begin{tabular}{c|c c c c c c c c c c}
edges & \multicolumn{10}{c}{number of vertices $n$} \\
$e$ & 12 & 13 & 14 & 15 & 16 & 17 & 18 & 19 & 20 & 21\\
\hline
6 & \textbf{1} &  &  &  &  &  &  &  &  &  \\
9 &  & \textbf{1} &  &  &  &  &  &  &  &  \\
12 &  &  & \textbf{1} &  &  &  &  &  &  &  \\
13 &  &  & 1 &  &  &  &  &  &  &  \\
15 &  &  &  & \textbf{1} &  &  &  &  &  &  \\
16 &  &  &  & 1 &  &  &  &  &  &  \\
17 &  &  &  & 1 &  &  &  &  &  &  \\
18 &  &  &  & 2 &  &  &  &  &  &  \\
20 &  &  &  &  & \textbf{2} &  &  &  &  &  \\
21 &  &  &  &  & $3_{(15)}$ &  &  &  &  &  \\
22 &  &  &  &  & $4_{(201)}$ &  &  &  &  &  \\
25 &  &  &  &  &  & \textbf{2} &  &  &  &  \\
26 &  &  &  &  &  & $7_{(30)}$ &  &  &  &  \\
27 &  &  &  &  &  & $6_{(642)}$ &  &  &  &  \\
28 &  &  &  &  &  & $1_{(13334)}$ &  &  &  &  \\
30 &  &  &  &  &  &  & \textbf{1} &  &  &  \\
31 &  &  &  &  &  &  & $7_{(15)}$ &  &  &  \\
32 &  &  &  &  &  &  & $14_{(382)}$ &  &  &  \\
33 &  &  &  &  &  &  & $4_{(8652)}$ &  &  &  \\
34 &  &  &  &  &  &  & $2_{(160573)}$ &  &  &  \\
37 &  &  &  &  &  &  &  & \textbf{11} &  &  \\
38 &  &  &  &  &  &  &  & $15_{(417)}$ &  &  \\
44 &  &  &  &  &  &  &  &  & \textbf{15} &  \\
51 &  &  &  &  &  &  &  &  &  & \textbf{4} \\
\end{tabular}
\caption{Number of $H_{13}$-patterned $(3,7)$-graphs.}
\end{table}
\end{center}

\begin{center}
\begin{table}
\tiny
\begin{tabular}{c|c c c c c c c c c c c}
edges & \multicolumn{11}{c}{number of vertices $n$} \\
$e$ & 14 & 15 & 16 & 17 & 18 & 19 & 20 & 21 & 22 & 23 & 24\\
\hline
7 & \textbf{1} &  &  &  &  &  &  &  &  &  &  \\
10 &  & \textbf{1} &  &  &  &  &  &  &  &  &  \\
13 &  &  & \textbf{1} &  &  &  &  &  &  &  &  \\
14 &  &  & 1 &  &  &  &  &  &  &  &  \\
16 &  &  &  & \textbf{1} &  &  &  &  &  &  &  \\
17 &  &  &  & 1 &  &  &  &  &  &  &  \\
18 &  &  &  & 1 &  &  &  &  &  &  &  \\
19 &  &  &  & 2 &  &  &  &  &  &  &  \\
20 &  &  &  &  & \textbf{1} &  &  &  &  &  &  \\
21 &  &  &  &  & 2 &  &  &  &  &  &  \\
22 &  &  &  &  & 3 &  &  &  &  &  &  \\
23 &  &  &  &  & 4 &  &  &  &  &  &  \\
25 &  &  &  &  &  & \textbf{2} &  &  &  &  &  \\
26 &  &  &  &  &  & 7 &  &  &  &  &  \\
27 &  &  &  &  &  & 7 &  &  &  &  &  \\
28 &  &  &  &  &  & 6 &  &  &  &  &  \\
29 &  &  &  &  &  & 1 &  &  &  &  &  \\
30 &  &  &  &  &  &  & \textbf{3} &  &  &  &  \\
31 &  &  &  &  &  &  & 13 &  &  &  &  \\
32 &  &  &  &  &  &  & 21 &  &  &  &  \\
33 &  &  &  &  &  &  & 14 &  &  &  &  \\
34 &  &  &  &  &  &  & 4 &  &  &  &  \\
35 &  &  &  &  &  &  & 2 & \textbf{1} &  &  &  \\
36 &  &  &  &  &  &  &  & 10 &  &  &  \\
37 &  &  &  &  &  &  &  & 48 &  &  &  \\
38 &  &  &  &  &  &  &  & 34 &  &  &  \\
39 &  &  &  &  &  &  &  & 15 &  &  &  \\
42 &  &  &  &  &  &  &  &  & \textbf{21} &  &  \\
43 &  &  &  &  &  &  &  &  & 100 &  &  \\
44 &  &  &  &  &  &  &  &  & 31 &  &  \\
45 &  &  &  &  &  &  &  &  & 15 &  &  \\
49 &  &  &  &  &  &  &  &  &  & \textbf{102} &  \\
50 &  &  &  &  &  &  &  &  &  & 51 &  \\
52 &  &  &  &  &  &  &  &  &  & 4 &  \\
56 &  &  &  &  &  &  &  &  &  &  & \textbf{51} \\
\end{tabular}
\caption{Number of $H_{13}$-patterned $(3,8)$-graphs.}
\end{table}
\end{center}

\begin{center}
\begin{table}
\tiny
\begin{tabular}{c|c c c c c c c c c c c c c}
edges & \multicolumn{13}{c}{number of vertices $n$} \\
$e$ & 16 & 17 & 18 & 19 & 20 & 21 & 22 & 23 & 24 & 25 & 26 & 27 & 28\\
\hline
8 & \textbf{1} &  &  &  &  &  &  &  &  &  &  &  &  \\
11 &  & \textbf{1} &  &  &  &  &  &  &  &  &  &  &  \\
14 &  &  & \textbf{1} &  &  &  &  &  &  &  &  &  &  \\
15 &  &  & 1 &  &  &  &  &  &  &  &  &  &  \\
17 &  &  &  & \textbf{1} &  &  &  &  &  &  &  &  &  \\
18 &  &  &  & 1 &  &  &  &  &  &  &  &  &  \\
19 &  &  &  & 1 &  &  &  &  &  &  &  &  &  \\
20 &  &  &  & 2 & \textbf{1} &  &  &  &  &  &  &  &  \\
21 &  &  &  &  & 1 &  &  &  &  &  &  &  &  \\
22 &  &  &  &  & 2 &  &  &  &  &  &  &  &  \\
23 &  &  &  &  & 3 &  &  &  &  &  &  &  &  \\
24 &  &  &  &  & 4 &  &  &  &  &  &  &  &  \\
25 &  &  &  &  &  & \textbf{2} &  &  &  &  &  &  &  \\
26 &  &  &  &  &  & 4 &  &  &  &  &  &  &  \\
27 &  &  &  &  &  & 7 &  &  &  &  &  &  &  \\
28 &  &  &  &  &  & 7 &  &  &  &  &  &  &  \\
29 &  &  &  &  &  & 6 &  &  &  &  &  &  &  \\
30 &  &  &  &  &  & 1 & \textbf{4} &  &  &  &  &  &  \\
31 &  &  &  &  &  &  & 14 &  &  &  &  &  &  \\
32 &  &  &  &  &  &  & 20 &  &  &  &  &  &  \\
33 &  &  &  &  &  &  & 21 &  &  &  &  &  &  \\
34 &  &  &  &  &  &  & 14 &  &  &  &  &  &  \\
35 &  &  &  &  &  &  & 4 & \textbf{4} &  &  &  &  &  \\
36 &  &  &  &  &  &  & 2 & 26 &  &  &  &  &  \\
37 &  &  &  &  &  &  &  & 62 &  &  &  &  &  \\
38 &  &  &  &  &  &  &  & 58 &  &  &  &  &  \\
39 &  &  &  &  &  &  &  & 34 &  &  &  &  &  \\
40 &  &  &  &  &  &  &  & 15 & \textbf{2} &  &  &  &  \\
41 &  &  &  &  &  &  &  &  & 13 &  &  &  &  \\
42 &  &  &  &  &  &  &  &  & 127 &  &  &  &  \\
43 &  &  &  &  &  &  &  &  & 184 &  &  &  &  \\
44 &  &  &  &  &  &  &  &  & 115 &  &  &  &  \\
45 &  &  &  &  &  &  &  &  & 31 &  &  &  &  \\
46 &  &  &  &  &  &  &  &  & 15 & \textbf{1} &  &  &  \\
47 &  &  &  &  &  &  &  &  &  & 35 &  &  &  \\
48 &  &  &  &  &  &  &  &  &  & 332 &  &  &  \\
49 &  &  &  &  &  &  &  &  &  & 412 &  &  &  \\
50 &  &  &  &  &  &  &  &  &  & 123 &  &  &  \\
51 &  &  &  &  &  &  &  &  &  & 51 &  &  &  \\
52 &  &  &  &  &  &  &  &  &  &  & \textbf{1} &  &  \\
53 &  &  &  &  &  &  &  &  &  & 4 &  &  &  \\
54 &  &  &  &  &  &  &  &  &  &  & 436 &  &  \\
55 &  &  &  &  &  &  &  &  &  &  & 705 &  &  \\
56 &  &  &  &  &  &  &  &  &  &  & 55 &  &  \\
57 &  &  &  &  &  &  &  &  &  &  & 51 &  &  \\
61 &  &  &  &  &  &  &  &  &  &  &  & \textbf{700} &  \\
62 &  &  &  &  &  &  &  &  &  &  &  & 95 &  \\
68 &  &  &  &  &  &  &  &  &  &  &  &  & \textbf{126} \\
\end{tabular}
\caption{Number of $H_{13}$-patterned $(3,9)$-graphs.}
\end{table}
\end{center}

\begin{center}
\begin{table}
\tiny
\begin{tabular}{c|c c c c c c c c c c c c c c}
edges & \multicolumn{14}{c}{number of vertices $n$} \\
$e$ & 18 & 19 & 20 & 21 & 22 & 23 & 24 & 25 & 26 & 27 & 28 & 29 & 30 & 31\\
\hline
9 & \textbf{1} &  &  &  &  &  &  &  &  &  &  &  &  &  \\
12 &  & \textbf{1} &  &  &  &  &  &  &  &  &  &  &  &  \\
15 &  &  & \textbf{1} &  &  &  &  &  &  &  &  &  &  &  \\
16 &  &  & 1 &  &  &  &  &  &  &  &  &  &  &  \\
18 &  &  &  & \textbf{1} &  &  &  &  &  &  &  &  &  &  \\
19 &  &  &  & 1 &  &  &  &  &  &  &  &  &  &  \\
20 &  &  &  & 1 &  &  &  &  &  &  &  &  &  &  \\
21 &  &  &  & 2 & \textbf{1} &  &  &  &  &  &  &  &  &  \\
22 &  &  &  &  & 1 &  &  &  &  &  &  &  &  &  \\
23 &  &  &  &  & 2 &  &  &  &  &  &  &  &  &  \\
24 &  &  &  &  & 3 &  &  &  &  &  &  &  &  &  \\
25 &  &  &  &  & 4 & \textbf{1} &  &  &  &  &  &  &  &  \\
26 &  &  &  &  &  & 2 &  &  &  &  &  &  &  &  \\
27 &  &  &  &  &  & 4 &  &  &  &  &  &  &  &  \\
28 &  &  &  &  &  & 7 &  &  &  &  &  &  &  &  \\
29 &  &  &  &  &  & 7 &  &  &  &  &  &  &  &  \\
30 &  &  &  &  &  & 6 & \textbf{3} &  &  &  &  &  &  &  \\
31 &  &  &  &  &  & 1 & 9 &  &  &  &  &  &  &  \\
32 &  &  &  &  &  &  & 14 &  &  &  &  &  &  &  \\
33 &  &  &  &  &  &  & 20 &  &  &  &  &  &  &  \\
34 &  &  &  &  &  &  & 21 &  &  &  &  &  &  &  \\
35 &  &  &  &  &  &  & 14 & \textbf{5} &  &  &  &  &  &  \\
36 &  &  &  &  &  &  & 4 & 27 &  &  &  &  &  &  \\
37 &  &  &  &  &  &  & 2 & 49 &  &  &  &  &  &  \\
38 &  &  &  &  &  &  &  & 62 &  &  &  &  &  &  \\
39 &  &  &  &  &  &  &  & 58 &  &  &  &  &  &  \\
40 &  &  &  &  &  &  &  & 34 & \textbf{5} &  &  &  &  &  \\
41 &  &  &  &  &  &  &  & 15 & 43 &  &  &  &  &  \\
42 &  &  &  &  &  &  &  &  & 166 &  &  &  &  &  \\
43 &  &  &  &  &  &  &  &  & 208 &  &  &  &  &  \\
44 &  &  &  &  &  &  &  &  & 184 &  &  &  &  &  \\
45 &  &  &  &  &  &  &  &  & 115 & \textbf{2} &  &  &  &  \\
46 &  &  &  &  &  &  &  &  & 31 & 19 &  &  &  &  \\
47 &  &  &  &  &  &  &  &  & 15 & 250 &  &  &  &  \\
48 &  &  &  &  &  &  &  &  &  & 734 &  &  &  &  \\
49 &  &  &  &  &  &  &  &  &  & 578 &  &  &  &  \\
50 &  &  &  &  &  &  &  &  &  & 412 &  &  &  &  \\
51 &  &  &  &  &  &  &  &  &  & 123 & \textbf{1} &  &  &  \\
52 &  &  &  &  &  &  &  &  &  & 51 & 63 &  &  &  \\
53 &  &  &  &  &  &  &  &  &  &  & 829 &  &  &  \\
54 &  &  &  &  &  &  &  &  &  & 4 & 2292 &  &  &  \\
55 &  &  &  &  &  &  &  &  &  &  & 1117 &  &  &  \\
56 &  &  &  &  &  &  &  &  &  &  & 705 &  &  &  \\
57 &  &  &  &  &  &  &  &  &  &  & 55 &  &  &  \\
58 &  &  &  &  &  &  &  &  &  &  & 51 & $4_{(5)}$ &  &  \\
59 &  &  &  &  &  &  &  &  &  &  &  & 1304 &  &  \\
60 &  &  &  &  &  &  &  &  &  &  &  & 4815 &  &  \\
61 &  &  &  &  &  &  &  &  &  &  &  & 1650 &  &  \\
62 &  &  &  &  &  &  &  &  &  &  &  & 700 &  &  \\
63 &  &  &  &  &  &  &  &  &  &  &  & 95 &  &  \\
66 &  &  &  &  &  &  &  &  &  &  &  &  & $5082_{(5084)}$ &  \\
67 &  &  &  &  &  &  &  &  &  &  &  &  & 2977 &  \\
69 &  &  &  &  &  &  &  &  &  &  &  &  & 126 &  \\
73 &  &  &  &  &  &  &  &  &  &  &  &  &  & \textbf{2657} \\
\end{tabular}
\caption{Number of $H_{13}$-patterned $(3,10)$-graphs.}
\end{table}
\end{center}

\begin{center}
\begin{table}
\tiny
\begin{tabular}{c|c c c c c c c c c c c c c c c c}
edges & \multicolumn{16}{c}{number of vertices $n$} \\
$e$ & 20 & 21 & 22 & 23 & 24 & 25 & 26 & 27 & 28 & 29 & 30 & 31 & 32 & 33 & 34 & 35\\
\hline
10 & \textbf{1} &  &  &  &  &  &  &  &  &  &  &  &  &  &  &  \\
13 &  & \textbf{1} &  &  &  &  &  &  &  &  &  &  &  &  &  &  \\
16 &  &  & \textbf{1} &  &  &  &  &  &  &  &  &  &  &  &  &  \\
17 &  &  & 1 &  &  &  &  &  &  &  &  &  &  &  &  &  \\
19 &  &  &  & \textbf{1} &  &  &  &  &  &  &  &  &  &  &  &  \\
20 &  &  &  & 1 &  &  &  &  &  &  &  &  &  &  &  &  \\
21 &  &  &  & 1 &  &  &  &  &  &  &  &  &  &  &  &  \\
22 &  &  &  & 2 & \textbf{1} &  &  &  &  &  &  &  &  &  &  &  \\
23 &  &  &  &  & 1 &  &  &  &  &  &  &  &  &  &  &  \\
24 &  &  &  &  & 2 &  &  &  &  &  &  &  &  &  &  &  \\
25 &  &  &  &  & 3 & \textbf{1} &  &  &  &  &  &  &  &  &  &  \\
26 &  &  &  &  & 4 & 1 &  &  &  &  &  &  &  &  &  &  \\
27 &  &  &  &  &  & 2 &  &  &  &  &  &  &  &  &  &  \\
28 &  &  &  &  &  & 4 &  &  &  &  &  &  &  &  &  &  \\
29 &  &  &  &  &  & 7 &  &  &  &  &  &  &  &  &  &  \\
30 &  &  &  &  &  & 7 & \textbf{2} &  &  &  &  &  &  &  &  &  \\
31 &  &  &  &  &  & 6 & 5 &  &  &  &  &  &  &  &  &  \\
32 &  &  &  &  &  & 1 & 9 &  &  &  &  &  &  &  &  &  \\
33 &  &  &  &  &  &  & 14 &  &  &  &  &  &  &  &  &  \\
34 &  &  &  &  &  &  & 20 &  &  &  &  &  &  &  &  &  \\
35 &  &  &  &  &  &  & 21 & \textbf{5} &  &  &  &  &  &  &  &  \\
36 &  &  &  &  &  &  & 14 & 18 &  &  &  &  &  &  &  &  \\
37 &  &  &  &  &  &  & 4 & 34 &  &  &  &  &  &  &  &  \\
38 &  &  &  &  &  &  & 2 & 49 &  &  &  &  &  &  &  &  \\
39 &  &  &  &  &  &  &  & 62 &  &  &  &  &  &  &  &  \\
40 &  &  &  &  &  &  &  & 58 & \textbf{8} &  &  &  &  &  &  &  \\
41 &  &  &  &  &  &  &  & 34 & 50 &  &  &  &  &  &  &  \\
42 &  &  &  &  &  &  &  & 15 & 133 &  &  &  &  &  &  &  \\
43 &  &  &  &  &  &  &  &  & 182 &  &  &  &  &  &  &  \\
44 &  &  &  &  &  &  &  &  & 208 &  &  &  &  &  &  &  \\
45 &  &  &  &  &  &  &  &  & 184 & \textbf{7} &  &  &  &  &  &  \\
46 &  &  &  &  &  &  &  &  & 115 & 68 &  &  &  &  &  &  \\
47 &  &  &  &  &  &  &  &  & 31 & 392 &  &  &  &  &  &  \\
48 &  &  &  &  &  &  &  &  & 15 & 694 &  &  &  &  &  &  \\
49 &  &  &  &  &  &  &  &  &  & 782 &  &  &  &  &  &  \\
50 &  &  &  &  &  &  &  &  &  & 578 & \textbf{3} &  &  &  &  &  \\
51 &  &  &  &  &  &  &  &  &  & 412 & 29 &  &  &  &  &  \\
52 &  &  &  &  &  &  &  &  &  & 123 & 469 &  &  &  &  &  \\
53 &  &  &  &  &  &  &  &  &  & 51 & 2199 &  &  &  &  &  \\
54 &  &  &  &  &  &  &  &  &  &  & 2795 &  &  &  &  &  \\
55 &  &  &  &  &  &  &  &  &  & 4 & 2424 &  &  &  &  &  \\
56 &  &  &  &  &  &  &  &  &  &  & 1117 & \textbf{1} &  &  &  &  \\
57 &  &  &  &  &  &  &  &  &  &  & 705 & 104 &  &  &  &  \\
58 &  &  &  &  &  &  &  &  &  &  & 55 & 1826 &  &  &  &  \\
59 &  &  &  &  &  &  &  &  &  &  & 51 & 8625 &  &  &  &  \\
60 &  &  &  &  &  &  &  &  &  &  &  & 8362 &  &  &  &  \\
61 &  &  &  &  &  &  &  &  &  &  &  & 5227 &  &  &  &  \\
62 &  &  &  &  &  &  &  &  &  &  &  & 1650 &  &  &  &  \\
63 &  &  &  &  &  &  &  &  &  &  &  & 700 & 11 &  &  &  \\
64 &  &  &  &  &  &  &  &  &  &  &  & 95 & 3403 &  &  &  \\
65 &  &  &  &  &  &  &  &  &  &  &  &  & 21553 &  &  &  \\
66 &  &  &  &  &  &  &  &  &  &  &  &  & 19748 &  &  &  \\
67 &  &  &  &  &  &  &  &  &  &  &  &  & 6330 &  &  &  \\
68 &  &  &  &  &  &  &  &  &  &  &  &  & 2977 &  &  &  \\
70 &  &  &  &  &  &  &  &  &  &  &  &  & 126 & 15 &  &  \\
71 &  &  &  &  &  &  &  &  &  &  &  &  &  & 26060 &  &  \\
72 &  &  &  &  &  &  &  &  &  &  &  &  &  & 37138 &  &  \\
73 &  &  &  &  &  &  &  &  &  &  &  &  &  & 3229 &  &  \\
74 &  &  &  &  &  &  &  &  &  &  &  &  &  & 2657 &  &  \\
77 &  &  &  &  &  &  &  &  &  &  &  &  &  &  & 4 &  \\
78 &  &  &  &  &  &  &  &  &  &  &  &  &  &  & 35173 &  \\
79 &  &  &  &  &  &  &  &  &  &  &  &  &  &  & 5314 &  \\
85 &  &  &  &  &  &  &  &  &  &  &  &  &  &  &  & 4307 \\
\end{tabular}
\caption{Number of $H_{13}$-patterned $(3,11)$-graphs.}
\end{table}
\end{center}

\begin{center}
\begin{table}
\tiny
\begin{tabular}{c|c c c c c c c c c c c c c c c c c}
edges & \multicolumn{17}{c}{number of vertices $n$} \\
$e$ & 22 & 23 & 24 & 25 & 26 & 27 & 28 & 29 & 30 & 31 & 32 & 33 & 34 & 35 & 36 & 37 & 38\\
\hline
11 & \textbf{1} &  &  &  &  &  &  &  &  &  &  &  &  &  &  &  &  \\
14 &  & \textbf{1} &  &  &  &  &  &  &  &  &  &  &  &  &  &  &  \\
17 &  &  & \textbf{1} &  &  &  &  &  &  &  &  &  &  &  &  &  &  \\
18 &  &  & 1 &  &  &  &  &  &  &  &  &  &  &  &  &  &  \\
20 &  &  &  & \textbf{1} &  &  &  &  &  &  &  &  &  &  &  &  &  \\
21 &  &  &  & 1 &  &  &  &  &  &  &  &  &  &  &  &  &  \\
22 &  &  &  & 1 &  &  &  &  &  &  &  &  &  &  &  &  &  \\
23 &  &  &  & 2 & \textbf{1} &  &  &  &  &  &  &  &  &  &  &  &  \\
24 &  &  &  &  & 1 &  &  &  &  &  &  &  &  &  &  &  &  \\
25 &  &  &  &  & 2 &  &  &  &  &  &  &  &  &  &  &  &  \\
26 &  &  &  &  & 3 & \textbf{1} &  &  &  &  &  &  &  &  &  &  &  \\
27 &  &  &  &  & 4 & 1 &  &  &  &  &  &  &  &  &  &  &  \\
28 &  &  &  &  &  & 2 &  &  &  &  &  &  &  &  &  &  &  \\
29 &  &  &  &  &  & 4 &  &  &  &  &  &  &  &  &  &  &  \\
30 &  &  &  &  &  & 7 & \textbf{1} &  &  &  &  &  &  &  &  &  &  \\
31 &  &  &  &  &  & 7 & 2 &  &  &  &  &  &  &  &  &  &  \\
32 &  &  &  &  &  & 6 & 5 &  &  &  &  &  &  &  &  &  &  \\
33 &  &  &  &  &  & 1 & 9 &  &  &  &  &  &  &  &  &  &  \\
34 &  &  &  &  &  &  & 14 &  &  &  &  &  &  &  &  &  &  \\
35 &  &  &  &  &  &  & 20 & \textbf{3} &  &  &  &  &  &  &  &  &  \\
36 &  &  &  &  &  &  & 21 & 10 &  &  &  &  &  &  &  &  &  \\
37 &  &  &  &  &  &  & 14 & 18 &  &  &  &  &  &  &  &  &  \\
38 &  &  &  &  &  &  & 4 & 34 &  &  &  &  &  &  &  &  &  \\
39 &  &  &  &  &  &  & 2 & 49 &  &  &  &  &  &  &  &  &  \\
40 &  &  &  &  &  &  &  & 62 & \textbf{7} &  &  &  &  &  &  &  &  \\
41 &  &  &  &  &  &  &  & 58 & 36 &  &  &  &  &  &  &  &  \\
42 &  &  &  &  &  &  &  & 34 & 77 &  &  &  &  &  &  &  &  \\
43 &  &  &  &  &  &  &  & 15 & 133 &  &  &  &  &  &  &  &  \\
44 &  &  &  &  &  &  &  &  & 182 &  &  &  &  &  &  &  &  \\
45 &  &  &  &  &  &  &  &  & 208 & \textbf{10} &  &  &  &  &  &  &  \\
46 &  &  &  &  &  &  &  &  & 184 & 86 &  &  &  &  &  &  &  \\
47 &  &  &  &  &  &  &  &  & 115 & 312 &  &  &  &  &  &  &  \\
48 &  &  &  &  &  &  &  &  & 31 & 520 &  &  &  &  &  &  &  \\
49 &  &  &  &  &  &  &  &  & 15 & 694 &  &  &  &  &  &  &  \\
50 &  &  &  &  &  &  &  &  &  & 782 & \textbf{9} &  &  &  &  &  &  \\
51 &  &  &  &  &  &  &  &  &  & 578 & 101 &  &  &  &  &  &  \\
52 &  &  &  &  &  &  &  &  &  & 412 & 813 &  &  &  &  &  &  \\
53 &  &  &  &  &  &  &  &  &  & 123 & 2156 &  &  &  &  &  &  \\
54 &  &  &  &  &  &  &  &  &  & 51 & 2786 &  &  &  &  &  &  \\
55 &  &  &  &  &  &  &  &  &  &  & 2795 & \textbf{3} &  &  &  &  &  \\
56 &  &  &  &  &  &  &  &  &  & 4 & 2424 & 42 &  &  &  &  &  \\
57 &  &  &  &  &  &  &  &  &  &  & 1117 & 798 &  &  &  &  &  \\
58 &  &  &  &  &  &  &  &  &  &  & 705 & 5504 &  &  &  &  &  \\
59 &  &  &  &  &  &  &  &  &  &  & 55 & 11359 &  &  &  &  &  \\
60 &  &  &  &  &  &  &  &  &  &  & 51 & 11153 &  &  &  &  &  \\
61 &  &  &  &  &  &  &  &  &  &  &  & 8362 & \textbf{1} &  &  &  &  \\
62 &  &  &  &  &  &  &  &  &  &  &  & 5227 & 162 &  &  &  &  \\
63 &  &  &  &  &  &  &  &  &  &  &  & 1650 & 3648 &  &  &  &  \\
64 &  &  &  &  &  &  &  &  &  &  &  & 700 & 25863 &  &  &  &  \\
65 &  &  &  &  &  &  &  &  &  &  &  & 95 & 45240 &  &  &  &  \\
66 &  &  &  &  &  &  &  &  &  &  &  &  & 31547 &  &  &  &  \\
67 &  &  &  &  &  &  &  &  &  &  &  &  & 19748 &  &  &  &  \\
68 &  &  &  &  &  &  &  &  &  &  &  &  & 6330 & 21 &  &  &  \\
69 &  &  &  &  &  &  &  &  &  &  &  &  & 2977 & 7765 &  &  &  \\
70 &  &  &  &  &  &  &  &  &  &  &  &  &  & 75617 &  &  &  \\
71 &  &  &  &  &  &  &  &  &  &  &  &  & 126 & 135323 &  &  &  \\
72 &  &  &  &  &  &  &  &  &  &  &  &  &  & 59825 &  &  &  \\
73 &  &  &  &  &  &  &  &  &  &  &  &  &  & 37138 &  &  &  \\
74 &  &  &  &  &  &  &  &  &  &  &  &  &  & 3229 &  &  &  \\
75 &  &  &  &  &  &  &  &  &  &  &  &  &  & 2657 & 102 &  &  \\
76 &  &  &  &  &  &  &  &  &  &  &  &  &  &  & 104844 &  &  \\
77 &  &  &  &  &  &  &  &  &  &  &  &  &  &  & 286371 &  &  \\
78 &  &  &  &  &  &  &  &  &  &  &  &  &  &  & 88162 &  &  \\
79 &  &  &  &  &  &  &  &  &  &  &  &  &  &  & 35173 &  &  \\
80 &  &  &  &  &  &  &  &  &  &  &  &  &  &  & 5314 &  &  \\
82 &  &  &  &  &  &  &  &  &  &  &  &  &  &  &  & 51 &  \\
83 &  &  &  &  &  &  &  &  &  &  &  &  &  &  &  & 299189 &  \\
84 &  &  &  &  &  &  &  &  &  &  &  &  &  &  &  & 153884 &  \\
86 &  &  &  &  &  &  &  &  &  &  &  &  &  &  &  & 4307 &  \\
90 &  &  &  &  &  &  &  &  &  &  &  &  &  &  &  &  & 133261 \\
\end{tabular}
\caption{Number of $H_{13}$-patterned $(3,12)$-graphs.}
\end{table}
\end{center}

\end{document}